\documentclass{article}

\usepackage{amssymb}

\newtheorem{Lemma}{Lemma}

\newtheorem{Corollary}[Lemma]{Corollary}

\newtheorem{Fact}[Lemma]{Fact}

\newcommand{\proof}{\noindent{\sc Proof}:\par\nobreak}
\newcommand{\qed}{\ \hfill$\square$\smallskip}

\newcommand{\Th}{\hbox{\rm Th}\,}

\newcommand{\<}{\langle}

\renewcommand{\>}{\rangle}

\newcommand{\Q}{\mathbb{Q}}

\title{Embeddings and chains of free groups} 

\author{Eric Jaligot, Azadeh Neman}
\date{May 27, 2008} 

\begin{document}

\maketitle 

\begin{abstract}
We build two nonabelian $CSA$-groups in which maximal abelian subgroups 
are conjugate and divisible, as the countable unions of increasing chains of $CSA$-groups 
and by keeping the constructions as free as possible in each case. 
\end{abstract}

For $n\geq 1$, let $F_{n}$ denote the free group on $n$ generators. We view 
all groups $G$ as first-order structures $\<G,\cdot,\phantom{g}^{-1},1\>$, where 
$\cdot$, $\phantom{g}^{-1}$, and $1$ denote respectively the multiplication, 
the inverse, and the identity of the group. 
The following striking results are proved in a series of papers of Sela culminating in 
\cite{Sela2007}. 

\begin{Fact}\label{FactTarskietStable}
{\bf \cite{MR2249582, MR2238944, MR2238945, Sela2007}}\
\begin{itemize}
\item[$(1)$]
For any $2\leq n \leq m$, the natural embedding $F_{n}\leq F_{m}$ is an elementary 
embedding. 
\item[$(2)$]
For any $n\geq 2$, the (common) complete theory $\Th(F_{n})$ is stable. 
\end{itemize}
\end{Fact}

We refer to \cite{Hodges(book)93} for model theory in general, and to 
\cite{Poizat(book)87} and \cite{Wagner(Book)97} for stability theory and 
in particular stable groups. 

Let $F$ denote the free group over countably many generators. 
Fact \ref{FactTarskietStable} has the following corollary. 

\begin{Corollary}
The natural embeddings 
$F_{2}\leq  \cdots F_{n} \leq \cdots \leq F$ are all elementary. In particular 
each $F_{n}$ is an elementary substructure of $F$, and $\Th(F)$ is stable. 
\end{Corollary}

A {\em $CSA$-group} is a group in which maximal abelian subgroups $A$ 
are {\em malnormal}, i.e., such that $A\cap A^{g}\neq 1$ implies that $g$ is in $A$ 
for any element $g$ of the ambient group. The class of $CSA$-groups contains 
all free groups and is studied from various points of view. We refer to 
\cite{HoucineJaligot04, MR2395048} for a model theoretic approach in 
combination of questions concerning particular groups 
\cite{Cherlin79, Jaligot01, MR2124721}, and to 
\cite{MR2381022} for an approach more related to computational aspects 
in limit groups. 

We prove the following lemma on embeddings of torsion-free $CSA$-groups in 
which maximal abelian subgroups are cyclic. 

\begin{Lemma}\label{LemmeEmbed}
Let $k\geq 2$, and let $G$ be a countable torsion-free $CSA$-group in which maximal 
abelian subgroups are cyclic. Let $a_{0}$ be any generator of a maximal abelian 
subgroup of $G$. Then $G$ embeds into a torsion-free $CSA$-group $H=\<F,r\>$ in 
which maximal abelian subgroups are cyclic, where $F$ is a free-group over countably 
many generators and $r^{k}=a_{0}$, and where maximal 
abelian subgroups of $G$ are $H$-conjugate. In particular any element of $G$ has 
a $k$-th root in $H$. 
\end{Lemma}
\proof
$G$ has countably many maximal abelian subgroups, and countably many conjugacy classes 
of such maximal abelian subgroups, which can be enumerated by $i<\omega$. 

For each such conjugacy class, fix a maximal (cyclic) abelian subgroup $A_{i}$, and 
inside $A_{i}$ fix a generator $a_{i}$. Notice that we can always take for $a_{0}$ 
any given generator of a maximal abelian subgroup of $G$. 

We define inductively on $i$ an increasing family of supergroups $G_{i}$ of $G$ as follows. 
\begin{itemize}
\item
$G_{0}=G$. 
\item
$G_{i+1}$ is the $HNN$-extension $\<G_{i},t_{i}~|~a_{i}^{t_{i}}=a_{i+1}\>$. 
\end{itemize} 
We note that each $G_{i+1}$ is generated by $G_{0}$ together with the elements 
$t_{0}$, ..., $t_{i}$. 

Let now 
$$G_{\omega}=\bigcup_{i<\omega} G_{i}.$$ 

It is clear that $G_{\omega}$ is generated by $G_{0}$, together with the elements $t_{i}$'s. 
By construction, it is also clear that in $G_{\omega}$ any two distinct elements 
$t_{i}$ and $t_{j}$ satisfy no relation by \cite[Britton's  Lemma]{LyndonSchupp77}. 
In particular they generate a free group on two generators $t_{i}$ and $t_{j}$. 

We also see that in $G_{\omega}$ any two maximal abelian subgroups of $G_{0}$ are 
conjugate, actually by the subgroup of $G_{\omega}$ generated by all the elements 
$t_{i}$. In particular in $G_{\omega}$ one has 
$$G_{0}\subseteq \<a_{0}\>^{\<(t_{i})_{i<\omega}~|~\>}$$
and $G_{\omega}$ is generated by $\<a_{0}\>$ and ${\<(t_{i})_{i<\omega}~|~\>}$. 

Consider now an abelian torsion-free cyclic supergroup $R$ of $A_{0}=\<a_{0}\>$ 
generated by an element $r$ such that $r^{k}=a_{0}$. As $k\geq 2$ by assumption, 
$r$ does not belong to the cyclic subgroup $\<a_{0}\>$ of $R$. 

Now one can form the free product of $R$ and $G_{\omega}$ with amalgamated subgroup 
$\<a_{0}\>$, say 
$$G_{\omega +1}=R \ast_{\<a_{0}\>}G_{\omega}.$$
As $r^{k}=a_{0}$ and $G_{\omega}$ is generated by $a_{0}$ and the $t_{i}$'s, one gets 
that $G_{\omega + 1}$ is generated by the $t_{i}$'s together with $r$. Hence 
$$G_{\omega+1}=\<r,(t_{i})_{i<\omega}\>$$
and the second set of generators freely generate the free group $F$. 

Now the natural embedding 
$$G\simeq G_{0}\leq G_{\omega+1}\simeq H$$
is the desired embedding. 

We note that the class of $CSA$-groups is inductive, as pointed out in 
\cite{HoucineJaligot04}. This follows indeed from the fact that the class 
is axiomatizable by universal axioms. As $H$ is a direct limit of $CSA$-groups, 
it is also a $CSA$-group. We note that maximal abelian subgroups of $H$ are also 
cyclic by the results of \cite{HoucineJaligot04}. 
\qed

\bigskip
By Lemma \ref{LemmeEmbed}, one can find an infinite sequence of embeddings 

$$^{1}G~\leq~^{2}G~\leq \cdots \leq~^{k-1}G~\leq~^{k}G~\leq \cdots $$
where $^{1}G$ is a nonabelian countable free group and 
such that for each $k\geq 2$ the embedding $^{k-1}G~\leq~^{k}G$ is as in 
Lemma \ref{LemmeEmbed}, i.e., such that maximal abelian 
subgroups of $^{k-1}G$ are conjugate in $^{k}G$ and each element in $^{k-1}G$ 
is $k$-divsible in $^{k}G$. We note that each group $^{k}G$ 
is a torsion-free $CSA$-group in which maximal abelian subgroups are cyclic. 

Consider now the group 
$$G=\bigcup_{k\geq 1} {^{k}G}$$
We see that $G$ is a $CSA$-group, as the class of $CSA$-groups is inductive. 
As maximal abelian subgroups coincide with centralizers 
of nontrivial elements in $CSA$-groups, one sees that maximal abelian subgroups of $G$ 
are conjugate by construction. Again the construction implies that each element $g$ of 
$G$ is $n$-divisible for each $n$, and as maximal abelian subgroups coincide with 
centralizers of nontrivial elements one concludes that maximal abelian subgroups 
are divisible. 

Of course $G$ is nonabelian since it contains the nonabelian free group $^{1}G$, 
and we note that it is also torsion-free. We obtain thus a non-abelian $CSA$-group 
in which maximal abelian subgroups are conjugate and divisible, as a countable union 
of torsion-free $CSA$-groups in which maximal abelian subgroups are cyclic. 

As the free group is stable, one may wonder about the stability of the group $G$ 
built above. Recall that the {\em stable} sets are the sets defined by a formula 
$\varphi(\overline{x},\overline{y})$ for which there exists an integer $n$, 
called the {\em ladder index} of $\varphi$, bounding uniformly the size of 
{\em $n$-ladders} of $\varphi$, that is the sets of tuples 
$$(\overline{x_{1}},\cdots , \overline{x_{1}}~;~\overline{y_{1}},\cdots ,\overline{y_{n}})$$ 
such that $\varphi(\overline{x_{i}},\overline{x_{j}})$ is true 
if and only if $i\leq j$. We refer to \cite{Hodges(book)93}. We 
recall also that by Ramsey's theorem boolean combinations of stable sets are stable 
\cite[0.2.10]{Wagner(Book)97}, and that the replacement of variables by 
parameters obviously does not affect the stability of definable sets. 

To check the stability of quantifier-free definable sets, the following lemma may 
be relevant for the type of groups considered here. 

\begin{Fact}\label{FactStableSetsQFIndLim}
Let $\varphi(\overline{x},\overline{y})$ be a quantifier-free formula in the 
langage of groups, that is a boolean combination of equations into the variables involved 
in the tuples $\overline{x}$ and $\overline{y}$. If a group $G$ is a union of an 
increasing family of subgroups $G_{i}$, 
and if $\varphi$ defines a stable set in each $G_{i}$ with a uniform bound 
on the ladder indices of $\varphi$ in each $G_{i}$, then 
$\varphi$ defines a stable set in $G$. 
\end{Fact}
\proof
By assumption, there exists a uniform bound on the ladder index 
of $\varphi$ in $G_{i}$, when $i$ varies, and thus we find an 
integer $n$ such that no subgroup $G_{i}$ can have a ladder of size $n$. 

Now one sees that $G$ cannot have a ladder of size $n$ also, as otherwise all 
elements of the tuples involved in such a ladder would belong to one of the subgroups 
$G_{i}$, a contradiction as $\varphi$ is quantifier-free. 

This proves that $\varphi$ defines a stable set in $G$. 
\qed

\bigskip
Another construction of a non-abelian $CSA$-group in which maximal abelian 
subgroups are conjugate and divisible and with a temptative to keep certain 
parts of the stability of the free group is as follows. It now consists of adding all roots 
simultaneously. 

Instead of starting from the free group $F$, let us start with $^{1}G=\Q \ast F$. 
Fix $a_{0}$ the element corresponding to the element $1$ of $\Q$ (in additive 
notation). Passing from $^{1}G$ to $^{2}G$ is now done as follows. Enumerate by 
$a_{1}$, $a_{2}$, ... etc, generators of maximal abelian cyclic subgroups of $^{1}G$, 
picking up exactly one maximal abelian cyclic subgroup in each conjugacy class of such 
subgroups. Now embed $\<a_{0}\>$ in a copy of $\Q$, in such a way that 
$a_{0}$ represents $1$ in $\Q$ (in additive notation), and form the free product 
of $^{1}G$ and this new copy of $\Q$, with amalgamated subgroup $\<a_{0}\>$. 
One gets then a new $CSA$-group, and one can conjugate $a_{0}$ to $a_{1}$ by 
forming an appropriate $HNN$-extension, with an element $t_{0}$. One repeats 
this process as in Lemma \ref{LemmeEmbed}, obtaining $CSA$-groups at each step 
by the results of \cite{HoucineJaligot04}. Calling $^{2}G$ the union, one gets a 
$CSA$-group with one conjugacy class of divisible maximal abelian subgroups, 
and generated by one such subgroup and a countable free group (generated by 
all the $t_{i}$'s added when forming the successive $HNN$-extensions). 

One can then build similarly an infinite sequence of embeddings 
$$^{1}G~\leq~^{2}G~\leq \cdots \leq~^{k-1}G~\leq~^{k}G~\leq \cdots $$
such that for each $k\geq 2$ the embedding $^{k-1}G~\leq~^{k}G$ is as in the 
process described above. In particular each $^{k}G$ is generated by a divisible abelian 
group isomorphic to $\Q$ and a free group $F$. 
Again the union of all these groups is a $CSA$-group 
in which maximal abelian subgroups are conjugate and divisible. 

\bigskip
{\it We thank Vincent Guirardel, Abderezak Ould Houcine, Zlil Sela, and Alina Vdovina 
who pointed out a mistake in the first version of the present work.}

\bibliographystyle{alpha}
\bibliography{biblio}



\end{document}